\documentclass[twoside,12pt,a4paper]{amsart} 
\usepackage{amsmath,amsfonts,amsthm,eurosym}
\usepackage[colorlinks=true,linkcolor=blue,citecolor=blue,filecolor=blue,urlcolor=blue,pageanchor=true,plainpages=false,linktocpage]{hyperref}
\usepackage[english]{babel}

\numberwithin{equation}{section}

\renewcommand{\epsilon}{\varepsilon}
\renewcommand{\phi}{\varphi}
\renewcommand{\rho}{\varrho}
\renewcommand{\theta}{\vartheta}

\begin{document}
\baselineskip3.6ex


\title{}

\author{}

\date{}

%


%

%

\makeatletter
\def\@makefnmark{}
\makeatother
\newcommand{\myfootnote}[2]{\footnote{\textbf{#1}: #2}}
 \footnote {}
{\bfseries\centerline{On PNDP-manifold}}
\\
\\
\centerline{A. Pigazzini, C. $\ddot{O}$zel, P. Linker and S. Jafari}
\\
\\
\centerline{Dedicated to Professor Maximilian Ganster on the occasion his retirement}
\\
\centerline{and Renata Albizetti 08/11/1925 - 06/01/2020}
\\
\\{\bfseries \centerline {Abstract}} 
\\ 
\\
We provide a possible way of constructing new kinds of manifolds which we will call Partially Negative Dimensional Product manifold (PNDP-manifold for short).
\\
In particular a PNDP-manifold is an Einstein warped product manifold of special kind, where the base-manifold $B$ is a Remannian (or pseudo-Riemannian) product-manifold $B=\Pi_{i=1}^{q'}B_i \times \Pi_{i=(q'+1)}^{\widetilde q} B_i$, with $\Pi_{i=(q'+1)}^{\widetilde q} B_i$ an Einstein-manifold, and the fiber-manifold $F$ is a derived-differential-manifold (i.e., $F$ is the form: smooth manifold ($\mathbb{R}^d$)+ obstruction bundle, so it can admit negative dimension).
\\
Since the dimension of a PNDP-manifold is not related with the usual geometric concept of dimension, from the speculative and applicative point of view, we try to define this relation using the concept of desuspension to identify the PNDP with another kind of "object", introducing a new kind of hidden dimensions.
\\
\\
\\
{\bfseries \centerline {1. Introduction and Preliminaries}}
\\

The concept of negative dimensional space is already used in linguistic statistics \cite{Maslov}. Also in supersymmetric theories in Quantum Field Theory, negative dimensional 
spaces are used \cite{Cvita}.


Let  $E \cong M \times F$ be a fiber bundle with base space $M$ and its fiber $F$.
We will discuss now a case, where the fiber has negative dimension. Note that the total dimension
of the fiber bundle is given by the relation
$\dim E = \dim M + \dim F$.
We will consider the case, where the base manifold has greater positive dimension than the negative
dimension of the fiber is  $\dim M > - \dim F$. In this case, the dimension of the
total fiber bundle is still positive. Since the base manifold is obtained by projection of the fiber
bundle along the fiber by projection operator $\pi_F$, we have
$\pi_F E = M$
i.e. the projection of the lower-dimensional fiber-bundle along the fiber yields the higher
dimensional base manifold space. Therefore, the projection operator  $\pi_F$ along the negative-dimensional fiber, is a suspension operator that raises the dimension of topological spaces.
\\
Let $\dim F = -d,d>0$. Then, the operator $\pi_F$ will be a $d$-fold suspension (a
single suspension extends a topological space by dimension one). By a slight redefinition of the
notion of a fiber bundle to the new mathematical structure, the inverse fiber bundle, we are enable to
turn fiber bundles with negative-dimensional fibers into inverse fiber bundles with positive-dimensional fibers. More precisely, we define the inverse fiber bundle as follows:
- An inverse fiber bundle $E^*$ is locally isomorphic to the space
\\ $M^* \times^* F^*$
with some base manifold $M^*$ and a fiber $F^*$.

- However, we have a suspension along the fiber 
$\Sigma_{F^*}$ such that $\Sigma_{F^*} E^* = M^*$ instead of a projection along a fiber in an inverse fiber bundle.
- Because of the use of suspension instead of projections, the inverse cartesian product 
$\times^*$ is not the ordinary cartesian product. It can be defined functorially as follows: in the
category of fiber bundle topology  $Fib$, the ordinary cartesian product is a morphism between
some topological spaces (the objects in this category). By defining an inversion functor 
$\mathcal{F}$ that maps the category $Fib$ to the category of inverse fiber bundles 
$InvFib$ by mapping projection operators (are also morphisms) to suspensions 
$\mathcal{F}(\pi_F) = \Sigma_{F^*}$. Objects remain the same, but on the fiber, the sign of the
dimension is changed (negative fibers in $Fib$ change to positive ones in  $InvFib$). Thus,
it holds $\mathcal{F}(\times) = \times^*$. We define this product to be compatible with
suspensions instead of projections. The functor also preserves commutative diagrams that are
frequently used in the definition of fiber bundles.
It is also easy to generalize the above discussion to the case of negative fibers with arbitrary negative
dimension. We can treat a fiber bundle with negative-dimensional fiber as an inverse fiber bundle,
where the fibers will be positive-dimensional. Cartesian products, however, must be changed to
inverse cartesian products. The base manifold in inverse fiber bundles has the highest dimension,
since it is obtained by a suspension. Imposing local isomorphisms between an inverse fiber bundle to
the inverse cartesian product of base manifold with fiber, we can see that the inverse cartesian
product reduces dimensions by some amount. It holds:
$\dim E^* = \dim(E^* \times^* F^*) = \dim M^* -
\dim F^*$.
\\
{\bfseries Definition 1:} We say that dimension of $X$ is $-1$ if suspension of $X$, $\Sigma X$ is diffeomorphic to a point. By induction and suspension operation we can define all negative dimensions. Since if $X$ and $Y$ are difeomorphic then their suspensions $\Sigma X$ and $\Sigma Y$ are diffeomorphic. Also a metric on $X$ is compatible with $\Sigma X$. It is also well-known that differentials and suspensions are compatible on manifolds. This means that the differential on the suspension space $\Sigma X$ is induced from the differential on $X$.

The operator $\times^*$ acts similar as the topological quotient, but in such a way that it is
compatible with respect to above defined functor. The inverse cartesian product we can also view as
a quotient operator obtained by using the $\mathcal{F}$-functor. We can alternatively write: ×
$\times^*= /_{\mathcal{F}}$.
\\
\\
{\bfseries Definition 2:} In general, given an n-dimensional space $X$, the suspension $\Sigma X$ has dimension $n+1$. Thus, the operation of suspension creates a way of moving up in dimension. The inverse operation $\Sigma^{-1}$, is called desuspension. Therefore, given an $n$-dimensional space $X$, the desuspension $\Sigma ^{-1}X$  has dimension $n-1$, (see \cite{Wolcott}).
\\
\\
\\
{\bfseries \centerline {1.1 PNDP-manifolds}}
\\

Einstein warped product manifolds are a mathematical object that is considered in current research. We will treat an Einstein warped product manifold in this section of the paper, where we will introduce also the concept of negative dimensions. As "mathematical tool", for the fiber-manifold, we will use derived geometry and about this we recall that if $A \rightarrow M$ and $B \rightarrow M$ are two transversal submanifolds of codimension $a$ and $b$ respectively, then their intersection $C$ is again a submanifold, of codimension $a+b$. Derived geometry explains how to remove the transversality condition and make sense out of a nontransversal intersection $C$ as a derived smooth manifold of codimension $a+b$. In particular, $dim C = dim M - a - b$, and the latter number can be negative. So the obtained dimensions are not related to the usual geometry concept of "dimension", but they are "virtual dimensions", and therefore from a speculative/applicative point of view we will try to relate the dimensions obtained from the PNDP-manifolds with desuspensions (special projections), interpreting this, as the negative dimensions of the fiber-manifold "hide" the positive dimensions of the base- manifold, making the PNDP-manifold appears as if desuspensions had been performed.
\\
\\
{\bfseries Definition 3:} A warped product manifold $(M, \bar{g})=(B,g)\times_f(F,\ddot{g})$ (where ($B, g$) is the base-manifold, ($F, \ddot{g}$) is the fiber-manifold), with $\bar{g}=g+f^2 \ddot{g}$,  is Einstein if only if:
\\
\\
\numberwithin{equation}{section}
{(1)}
$\bar{Ric}=\lambda \bar{g} \Longleftrightarrow\begin{cases} 
 Ric- \frac{d}{f}\nabla^2 f= \lambda g  \\  \ddot{Ric}=\mu \ddot{g} \\ f \Delta f+(d-1) |\nabla f|^2 + \lambda f^2 =\mu
\end{cases}$
\\
\\
where $\lambda$ and $\mu$ are constants, $d$ is the dimension of $F$, $\nabla ^2f$, $ \Delta f$ and $\nabla f$ are, 
\\
respectively, the Hessian, the Laplacian and the gradient of $f$ for $g$, with $f:(B) \rightarrow \mathbb{R}^+$ a smooth positive function.
\\
\\
Contracting first equation of (1) we get: 
\\
\\
\numberwithin{equation}{section}
{(2)}
$R_Bf^2-f \Delta fd=n f^2 \lambda$ 
\\
where $n$ and $R_B$ is the dimension and the scalar curvature of $B$ respectively, and from third equation, considering $d \neq 0$ and $d \neq 1$, we have:
\\
\\
\numberwithin{equation}{section}
{(3)}
$f\Delta fd+d(d-1)|\nabla f|^2+\lambda f^2d=\mu d$
\\
Now from (2) and (3) we obtain:
\\
\\
\numberwithin{equation}{section}
{(4)}
$|\nabla f|^2+[\frac{\lambda (d-n)+R_B}{d(d-1)}]f^2=\frac{\mu}{(d-1)}$.
\\
\\
{\bfseries Definition 4:}  We called PNDP-manifold a warped product manifold \\ $(M, \bar{g})=(B, g) \times_f (F, \ddot{g})$ that satisfies (1), where the base-manifold $(B, g)$ is a Riemannian (or pseudo-Riemannian) product-manifold $B=\Pi_{i=1}^qB_i=B_1 \times B_2 \times ...$ with $g=\Sigma g_i$, which we can write as $B=B'\times \widetilde B$ where $\widetilde B$ is an Einstein manifold $\Pi_{i=(q'+1)}^{\widetilde q} B_i$ (i.e., $\widetilde Ric=\lambda \widetilde g$ where $\lambda$ is the same for (1) and $\widetilde g$ is the metric for $\widetilde B$), with $dim\widetilde B=\widetilde n$ and $B'=\Pi_{i=1}^{q'}B_i$ with $dimB'=n'$, so $dimB=n=n'+\widetilde n$. The warping function $f:B \rightarrow \mathbb{R}^+$ is $f(x,y)=f'(x)+ \widetilde f(y)$  (where each is a function on its individual manifold, i.e., $f':B'\rightarrow \mathbb{R}^+$ and $\widetilde f:\widetilde B \rightarrow \mathbb{R}^+$) and can also be a constant function. The fiber-manifold $(F, \ddot{g})$ is a derived differential Riemann-flat manifold with negative integer dimensions $m$, where with \textit{derived differential manifold} we consider the derived manifolds as \textit{smooth Riemannian flat manifolds by adding a vector bundle of obstructions}. In particular we consider, for $F$, only the spaceforms $\mathbb{R}^d$, with orthogonal Cartesian coordinates such that $g_{ij}=-\delta_{ij}$, by adding a vector bundle of obstructions, $E \rightarrow \mathbb{R}^d$, (so $F:=\mathbb{R}^d + E$), with dimension $m=d - rank(E)$, where $rank(E) > d$ (for issues related to derived-geometry or obstruction bundle see \cite{Joyce}), and more specifically we consider $F$ such that its dimension is $m=-d$. Since $F=\mathbb{R}^d+E$, and on $E$ (obstruction bundle) the Riemannian geometry does not work, we consider and define that every operation of Riemannian geometry done and defined on the underlying $\mathbb{R}^d$, is considered done and defined on $F$ (i.e., for example, if we calculate the Ricci curvature of $\mathbb{R}^d$, which is obviously zero, then we will say that the Ricci curvature of $F$ is zero, this because it will be the definition of Ricci curvature for $F$). Moreover if we want the quantity $n-d>0$ then we must set $n'=d=-m$, and in the special case where $n-d>0$ with also $B'$ an Einstein-manifold with the same Einstein-$\lambda$, then $B'_i=\widetilde B_i$, or $B'=\widetilde B$.
\\
\\
{\bfseries Important Note 1:} Since the usual Riemannian geometry works for the underlying Riemannian manifold as for "ordinary" manifolds, we can, therefore, consider ourselves to be able to work with the (pseudo-)Riemannian geometry on our derived fibers-manifold $F$, defining all Riemannian geometry operations on $\mathbb{R}^d$ as Riemannian geometry operations made on $F$, paying attention only to the dimension; In particular, as mentioned above, we work on $F$ by actually working on the underlying $\mathbb{R}^d$, (in the sense that, since $F:= \mathbb{R}^d + E$, the calculation results obtained on $\mathbb{R}^d$ are considered to be obtained and defined on $F$), but taking into account that $F$ has negative dimensions, in fact on the obstruction bundle we cannot put the metric, then we do not work on the obstruction bundle with the Riemannian geometry, but it gives negative dimensions on $F$. Obviously for what has been said, the tangent space and the vector fields are those of $\mathbb{R}^d$. The scalar product with two arbitrary vector fields ${^F}g\langle V, W \rangle$ is define on $F$ as: $g_{ij}v^iw^j=-\delta_{ij} v^i w^j= -(v^iw^i)$. 
\\
\\
We want to underline that the analysis does not differ from the usual Einstein warped product manifold analysis, and that the Riemannian curvature tensor and the Ricci curvature tensor of the product Riemannian manifold can be written respectively as the sum of the Riemannian curvature tensor and the Ricci curvature tensor of each Riemannian manifold (see \cite{Atceken}). This could be also a special kind of Einstein Sequential warped product manifold $(M_1 \times_h M_2) \times_{\bar{h}}M_3$, (see \cite{Chand Dea}), where $h=1$, $M_2$ is an Einstein-manifold and $M_3$ is a derived-smooth-manifold with negative dimensions, that said:
\\
\\
{\bfseries Proposition:} If we write the B-product as $B=B' \times \widetilde B$, where:
\\
i) $Ric'$ is the Ricci tensor of $B'$,  
\\
ii) $\widetilde Ric$ is the Ricci tensor of $\widetilde B$,
\\
iii) $g'$ is the metric tensor referred to $B'$, 
\\
iv) $\widetilde g$ is the metric tensor referred to $\widetilde B$,
\\
v) $f=f'+ \widetilde f$, is the smooth warping function $f:B \rightarrow \mathbb{R}^+$, where $f':B' \rightarrow \mathbb{R}^+$ and $\widetilde f: \widetilde B \rightarrow \mathbb{R}^+$,
\\
vi) $\nabla^2 f=\tau'^*\nabla'^2f'+\widetilde \tau^* \widetilde \nabla^2 \widetilde f$ is the Hessian referred on its individual metric, but since $\widetilde B$ is Einstein, we have $\widetilde \tau^* \widetilde \nabla^2 \widetilde f=0$, (where $\tau'^*$ is the pullback), 
\\
vii) $\nabla f$ is the gradient (then $|\nabla f|^2= |\nabla' f'|^2 + |\widetilde \nabla \widetilde f|^2$), and
\\
viii) $\Delta f=\Delta' f'+ \widetilde \Delta \widetilde f$ is the Laplacian, (since $\widetilde \tau^* \widetilde \nabla^2 \widetilde f=0$, then $\widetilde \Delta \widetilde f=0$), so the Ricci curvature tensor will be:
\\
\\
\numberwithin{equation}{section}
{(1****)}
$\begin{cases} 
\bar{Ric}(\Sigma_{i=1}^{q'}X_i, \Sigma_{i=1}^{q'}Y_i) = Ric'(\Sigma_{i=1}^{q'}X_i, \Sigma_{i=1}^{q'}Y_i) - \frac{d}{f}\nabla'^2 f'(\Sigma_{i=1}^{q'}X_i, \Sigma_{i=1}^{q'}Y_i)\\  \bar{Ric}( \Sigma_{i=(q'+1)}^{\widetilde q}X_i, \Sigma_{i=(q'+1)}^{\widetilde q}Y_i) = \widetilde Ric( \Sigma_{i=(q'+1)}^{\widetilde q}X_i, \Sigma_{i=(q'+1)}^{\widetilde q}Y_i)  \\ \bar {Ric}( U, V)=\ddot{Ric}( U, V) -f^2 \ddot{g}( U, V)f^*\\  \bar {Ric}(\Sigma_{i=1}^{q'}X_i, \Sigma_{i=(q'+1)}^{\widetilde q}X_i)=0 \\ \bar {Ric}(\Sigma_{i=1}^{q}X_i, U)=0, 
\end{cases}$
\\
where  $f^*= \frac{\Delta' f'}{f}+(d-1) \frac{|\nabla f|^2}{ f^2}$, $\Sigma_{i=1}^{q'}X_i=X'$, $\Sigma_{i=1}^{q'}Y_i=Y'$,  $\Sigma_{i=(q'+1)}^{\widetilde q}X_i=\widetilde X$ and $\Sigma_{i=(q'+1)}^{\widetilde q}Y_i=\widetilde Y$.
\\
\\
{\bfseries Theorem:} A warped product manifold with derived differential fiber-manifold \\ $F:=\mathbb{R}^d+E$, and $dimF$ a negative integer, is a  PNDP-manifold, as defined in \\ \textit{Definition 4}, if and only if:
\\
\\
\numberwithin{equation}{section}
{(1***)}
$\bar{Ric}=\lambda \bar{g} \Longleftrightarrow\begin{cases} 
 Ric'- \frac{d}{f}\tau'^*\nabla'^2 f'= \lambda g'\\ \widetilde \tau^* \widetilde \nabla^2 \widetilde f=0
\\
\widetilde Ric = \lambda \widetilde g \\ \ddot{Ric}=0 \\ f \Delta' f'+(d-1) |\nabla f|^2 + \lambda f^2 =0,
\end{cases}$
\\
(since $Ric$ is the Ricci curvature of $B$, then $Ric=Ric'+\widetilde Ric=\lambda(g'+\widetilde g)+ \frac{d}{f}\tau'^*\nabla'^2 f'$).
\\
\\
Therefore (2) and (3), for $n-d=0$ and $n-d<0$, become:
\\
\\
\numberwithin{equation}{section}
{(1**)}
$\bar{R}=\lambda \bar n \Longleftrightarrow\begin{cases} 
R'f- \Delta' f'd=n' f \lambda \\ \widetilde \Delta \widetilde f=0 \\ \widetilde R = \lambda \widetilde n \\ \ddot{Ric}=0 \\ f \Delta' f'+(d-1) |\nabla f|^2 + \lambda f^2 =0.
\end{cases}$
\\
\\
where $n'$ and $R'$ are the dimension and the scalar curvature of $B'$ respectively, while for $n-d>0$, we must set $d=n'$. We have
\\
\numberwithin{equation}{section}
{(1*)}
$\bar{R}=\lambda \bar n \Longleftrightarrow\begin{cases} 
R'f- \Delta' f'n'=n' f \lambda \\ \widetilde \Delta \widetilde f=0 \\ \widetilde R = \lambda \widetilde n \\ \ddot{Ric}=0 \\ f \Delta' f'+(n'-1) |\nabla f|^2 + \lambda f^2 =0.
\end{cases}$
\\
\\
\\
\textit{Proof.} We applied the condition that the warped product manifold of system (1****) is Einstein. $\square$
\\
\\
Recapitulating, we used the Derived-geometry to define the fiber-manifold, which explains how to remove the transversality condition and make sense of a non-transversal intersection and therefore admit the presence of negative dimensions. Since we consider the fiber-manifold $F$ as $\mathbb{R}^d$ adding a vector bundle of obstructions (then we work with the Riemannian geometry on $F$ working only on the underlying $\mathbb{R}^d$, such that all the Riemannian geometry operations made on $\mathbb{R}^d$ are considered as defined on $F$), then the classical construction for the warped product manifold (see \cite{O'Neill}) is the same; we consider the vertical vector fields $U$, $V$, as lift of vector fields of $F$ the development of the formulas remains the same; $X, Y$ are lift from $B$, they are horizontal and so constant on fibers, then for example $V[X, Y]=0$, and we continue to consider the inner product between a vector field on $B$ with one on $F$ as zero (i.e., $\langle X, V\rangle =0$). We obtain an Einstein warped product manifold  $(M, \bar{g})=(B, g) \times_f (\mathbb{R}^d, \ddot{g})$, because from the geometric point of view the fiber-manifold is $\mathbb{R}^d$, but taking into account that the manifold obtained has the fiber-manifold with negative dimensions, then the outcome manifold is $(n+(d-rank(E)))$-dimensional manifold.
\\
Another important observation is that we consider $\mathbb{R}^d$ because, from a speculative point of view, negative dimensional space is used only to make the PNDP-manifold dimension smaller than the base-manifold dimension. Therefore for this purpose we consider it in the simplest possible form.
\\
\\
{\bfseries PNDP-metric:} Referring to a PNDP-manifolds, with negative dimensional fiber, and for not confusing its metric with the metrics of a "classic" Einstein warped product manifold, we denote the Riemannian or pseudo-Riemannian metric of the fiber-manifold with the following notation to indicate that $F$ has negative dimensions: \\ $\bar{g}=g-f^2(\Sigma^n_{i=1}(d\psi^i)^2)_{(m)}$, where $g$ is the metric of the base-manifold $B$, and $m$ is the negative dimension of $F$. Then the general metric form of a PNDP-manifold is:
\\
$\bar g=g-f^2(\Sigma^n_{i=1}(d\psi^i)^2)_{(m)}=(g' + \widetilde g)+(f'+ \widetilde f)^2(\Sigma^n_{i=1}(d\psi^i)^2)_{(m)}$.
\\
\textit{Example 1 - A type of flat PNDP-manifold:} 
\\
The manifold $(\mathbb{R} \times \mathbb{R} \times \mathbb{R}\times \mathbb{R}) \times[(\mathbb{R}^2+E)]$ (with $rank(E)=4$) is a $(4-2)$-PNDP-manifold Ricci-flat. In fact it satisfies (1*) for constant $f=1$ ($f'$ and $\widetilde f$ both constants), we have $n+m=n-d=4-2=2>0$. Thus we have to consider $dimB'=dimF$. In this special case the B-product is composed of only Einstein-manifolds where $B'=\widetilde B$, that is $\mathbb{R}^2$, and its metric is: $ds^2=dt^2+dx^2+dy^2+dz^2-(du^2+dv^2)_{(-2)}$.
\\
\\
\textit{Example 2 - A type of PNDP-point-like manifold:} For this purpose, we consider fiber-manifold $F=(\mathbb{R}^3+E)$ (with $rank(E)=6$), 1-dimensional $B'$-manifold and $\lambda=0$. Like in the previous case, if for example we choose $\widetilde B$ as $\mathbb{R}^2$, then our PNDP-manifold will be: 
\\
$M=(\mathbb{R}^2 \times \mathbb{R})_{f}(\mathbb{R}+E)$, and with the following metric: 
\\
$ds^2=dx^2+dy^2+dz^2 + h(y)^2(du^2+dv^2+dw^2)_{(-3)}$.
\\
\\
\textit{Example 3 - A type of PNDP-manifold with negative "virtual" dimension:} Following what was done in the previous examples, keeping the same setting,  and considering, for example, $dim(F)=-4$, we obtain that space "will emerge" as a desuspension of a point, i.e., negative dimension. 
\\
\\
{\bfseries \centerline {2. Possible speculative interpretation of the types of PNDP manifolds}}
\\
\\
We assume that negative quantities exist in nature, here we consider the possible existence of negative spatial dimensions, but the fact remains that the PNDP-manifolds could be used to describe anything in nature that has to do with negative quantities.
\\
\\
From the speculative point of view, in this session we consider 2 types of PNDP-manifolds and try to relate their "virtual dimensions" to the usual concept of dimensions in Riemannian geometry, and this interpreting that, through the product, the negative virtual dimensions of the fiber-manifold "hide" the positive dimensions of the base-manifold, making the PNDP-manifold appears as if  desuspensions had been performed.
\\
In practice the derived-manifold $F$ has negative dimensions, which are far from the usual geometric concept of "dimension". Being that with the "warped product" we must consider both the dimension of $F$ (virtual) and the dimension of B, (which is a Riemannian or pseudo-Riemannian manifold and therefore its dimension is in perfect relationship with the usual concept of "dimension"), the obtained virtual dimension of our PNDP-manifold is equivalent to the dimension of the base-manifold $B$ with the dimension of the underlying manifolds $\mathbb{R}^d$ subtracted, and we get that the result will be smaller than the dimension of $B$. From the speculative point of view, we want to consider the dimensions obtained for the PNDP-manifold interpreting that when the dimensions of $B$ come into "contact" with the dimensions of $F$, these latter interact by "hiding" the dimensions of $B$. This means that the PNPD-manifold, will "visually appears" as if a kind of desuspensions had been performed on $B$, then with extra hidden dimensions which we can not see. 
\\
The same thing can be assumed for $F$; it "will appear" only as desuspensions of the point, with extra hidden dimensions which we can not see. Below we have consider two kinds of PNDP-manifolds:
\\
\\
{\bfseries Type I)} the $(n, -n)$-PNDP manifold that has overall, zero-dimension ($dim M = dim B + dim F = n + (-n) = 0$).
Then the result may be interpreted as an "invisible" manifold but made up of two manifolds with $n$ and $-n$ dimensions, respectively. Then we try to consider it as a kind of "point-like manifold" (zero-dimension), but with "hidden" dimensions, and
\\
\\
{\bfseries Type II)} the $(n,-d)$-PNDP manifold, where $n$ (the base-manifold dimension) is different from $d$ such that $dim=n-d>0$. The particular feature of this manifold is that it appears as another Einstein-manifold with $(n-d)$ dimension. 
\\
\\
In (Type I) we try to consider $dim=0$ referring to a point-like manifold. It would be thought that the points of $B$ and of the underlying $\mathbb{R}^d$ of $F$ appear as if they were "projected" on the PNDP-manifold, which in this case degenerates into a point (dim = 0), and we observe the manifold as a point-like; in effect this perception is given by our interpretation that negative dimensions "hide" positive dimensions and show us the manifold as if it were a point, i.e. the PNDP-manifold is "hidden" except for a point.
\\
Also in (Type II) case, since the dimension of the PNDP-manifold is calculated by the dimension of the base-manifold with the dimension of the fiber-manifold subtracted, we try to interpret it as if the negative dimensions "hide" the positive dimensions of the base-manifold. Since the total dimension will be smaller than the dimension of B, but still positive,we consider that the PNDP-manifold appears as a projection of $B$ with $dim = n-d$.
\\
Therefore our projection interpret this operation by considering that each negative dimension can act equally on each positive dimension, so that the final dimension is the difference between positive and negative dimensions. 
\\
In our speculative interpretation we consider $(n,-d)$-PNDP manifold like an Einstein-manifold projected from the product-manifold $B$, which has hidden dimensions, but we will see this part better in the next section.
\\
Concluding this section, we can say that our interpretation makes Type I and II "appear" as point or classic Einstein manifold, different from what they actually are.
\\
\\
\\
{\bfseries \centerline {3. Desuspension as Projection}}
\\
\\
This section helps to better define how the desuspensions work.
\\
As we have mentioned in the previous paragraph, we consider the relation between PNDP-manifold dimension and the usual geometric concept of "dimension", as desuspensions interpreted as projections, and for this, since each negative dimension can act equally on each positive dimension, in order not to be faced with the choice on which projection of B the result can be, we must "orient" the result of the projections to what reflects certain characteristics and to do this we must consider $B$ as a Riemannian base-manifold expressible as a Riemannian product-manifold $B = B '\times \widetilde B$ with $g_B = g_ {B'} + g _ {\widetilde B}$, such that $\widetilde B$ is an Einstein manifold with the same dimension and the same constant $\lambda$ of the PNDP-manifold, the latter will be the result of the projection we are going to consider, so:
\\
\\
$B= \Pi_{i=1}^{q'} B_i \times \Pi_{i=(q'+1)}^{\widetilde q}B_i$ and $F=\mathbb{R}^d + E$, then PNDP$=B \times_f F=[\Pi_{i \in I}B_i \times_f F]$, with $dimB -dim\mathbb{R}^d=dim$PNDP$=dim(\pi_{(n-d)}:$ PNDP$)=(n-d)$.
\\
- If $(n-d)>0$ (i.e. system solutions (1*)) we have the projection:
\\ 
$\pi_{(n-d)}:$PNDP$\rightarrow (\Pi_{i=(q'+1)}^{\widetilde q}B_i)=\widetilde B$, 
\\
- if $(n-d)=0$, (i.e., system solutions (1**)), we have the projection:
\\ 
$\pi_{0}:$PNDP$\rightarrow $ point-like manifold, and
\\
- if $(n-d)<0$, (i.e., system solutions (1**)), we have the projection:
\\
$\pi_{(n-d)}:$PNDP$\rightarrow \Sigma^{n-d}(p)$, with $\Sigma^{n-d}(p)$, we mean the $(n-d)$-th desuspension of point.
\\
\\
\\
\textit{Example 4:} Let $(B'_1 \times B'_2 \times B'_3 \times \widetilde B_4 \times \widetilde B_5 \times\widetilde B_6) \times_fF$ be a $(6-3)$-PNDP-manifold with $f$ non-constant, and since $n+m=n-d=6-3=3>0$, then $dimB'=dimF$.
\\
So our PNDP-manifold is such that: $\widetilde B_4 \times \widetilde B_5 \times \widetilde B_6$ will be an Einstein-manifold, i.e., $Ric_{(\widetilde B_4 \times \widetilde B_5 \times \widetilde B_6)}=\lambda (g_{\widetilde B_4}+g_{\widetilde B_5}+ g_{\widetilde B_6})$, and since $n-d=3$ we have the desuspension/projection: $\pi_3:$$(6-3)$-PNDP$\rightarrow (\Pi_{i=4}^6 B_i)$, by identifying the $(6-3)$-PNDP manifold with the Einstein-manifold $\widetilde B_4 \times \widetilde B_5 \times \widetilde B_6$.
\\
\\
\textit{Example 5:} Let $(B'_1 \times B'_2 \times \widetilde B_3 \times \widetilde B_4) \times_fF$ be a $(8-4)$-PNDP-manifold with $f$ non-constant, and also in this case, since $n+m=n-d=6-4=4>0$, $dimB'=dimF$.
\\
Then our PNDP-manifold is such that: $\widetilde B_3 \times \widetilde B_4$ will be Einstein, i.e., \\ $Ric_{(\widetilde B_3 \times \widetilde B_4)}=\lambda (g_{\widetilde B_3}+g_{\widetilde B_4})$, then we will have: $\pi_4:(8-4)$-PNDP $\rightarrow (\widetilde B_3 \times \widetilde B_4)$. Hence we identify the $(8-4)$-PNDP manifold with the Einstein-manifold $\widetilde B_3 \times \widetilde B_4$ . 
\\
\\
\textit{Example 6:} If we consider the manifold of the \textit{Example 1}, where we have the special case $B'=\widetilde B$, that is $\mathbb{R}^2$, the desuspension/projection will be: 
\\
$\pi_2:(\mathbb{R} \times \mathbb{R} \times \mathbb{R}\times \mathbb{R}) \times(\mathbb{R}^2+E)\rightarrow \mathbb{R}^2$, i.e., we identify the $(4-2)$-PNDP manifold with $\mathbb{R}^2$.
\\
\\
\textit{Example 7:}  If we consider the manifold of the \textit{Example 2}, the desuspension/projection will be:
$\pi_{0}:(\mathbb{R}^2 \times \mathbb{R}) \times_{f}(\mathbb{R}^3+E)\rightarrow$ point-like manifold (zero-dimension).
\\
\\
\textit{Example 8:}  If we consider the manifold of the \textit{Example 3}, the desuspension/projection will be:
$\pi_{-1}:(\mathbb{R}^2 \times \mathbb{R}) \times_{f}(\mathbb{R}^4+E)\rightarrow \Sigma^{-1}(p)$, i.e., a desuspension of a point.
\\
\\
{\bfseries Important Note 2:} In summary, from a speculative point of view, we hypothesize that in nature there are hidden dimensions that can be interpreted as the interaction of negative dimensions, (which we have introduced with manifolds generated by non-transversal intersections, derived-geometry), with positive dimensions, which by closing together, are not perceived, therefore hidden. To show this effect, we used special projections such as desuspensions, which show only a part of the manifold, its virtual dimension, while the remaining part is hidden. In this speculative context the PNDP is identified by $(\pi_{(n-d)}, \lambda, (n,m), g)$.
\\
\\
\\
\\
{\bfseries \centerline {4. Possible speculative application about the $(n, -n)$-PNDP-manifolds}}
\\
{\bfseries \centerline {and $(n, -q)$-PNDP-manifolds}}
\\
\\
In this section we present possible applications in light of the interpretations made above.
\\
Since $m=d - rank E$, where $rank E > d$ (such that $m=-d$), here we consider $d=n$, then $m=-n$, where $n$ is the dimension of the base-manifold. Referring to a manifold with $n$ negative dimensions, we denote the pseudo-Riemannian metric of the fiber-manifold with the following notation to indicate that $F$ has negative dimensions: $\ddot{g}=- (\Sigma^n_{i=1}(d\psi^i)^2)_{(-n)}$, where $_{(-n)}$ is the dimension of the fiber-manifold $F$. In this way we don't confuse it with the metrics of a "classic" product manifold.
\\
Then the $(n, -n)$-PNDP metric has the form: $\bar{g}=g-f^2(\Sigma^n_{i=1}(d\psi^i)^2)_{(-n)}$, where $g$ is the metric of the base-manifold $B$.
\\
First of all we recall and highlight that the purpose of the PNDP-manifolds is precisely to present the point-like manifolds from a mathematical point of view, and introduce a type of manifold with a new kind of hidden dimensions.
\\
In \cite{Capozziello}, Capozziello et al. introduced the concept of the "point-like manifold" building superconductors with graphene. In particular they argued that superconductor graphene can be produced by molecules organized in point-like structures where sheets are constituted by $(N+1)$-dimensional manifold. Particles like electrons, photons and ‘‘effective gravitons’’ are string modes moving on this manifold. In fact, according to string theory, bosonic and fermionic fields like electrons, photons and gravitons are
particular ‘‘states’’ or ‘‘modes’’ of strings. In their important work, they show that at the beginning, there are point-like polygonal manifolds (with zero spatial dimension) in space which strings are attached to them, where all interactions between strings on one manifold are the same and are concentrated on one point which the manifold is located on it. They also show that by joining these manifolds, 1-dimensional polygonal manifolds are emerged on which gauge fields and gravitons live and so, these manifolds glued to each other build higher dimensional polygonal manifolds with various orders of gauge fields and curvatures.
\\
As in \cite{Capozziello}, the authors in \cite{Sepehri}, propose a version of Moffat's Modified Gravity (MOG) without anomaly, where they showed the same "configuration" for space-time where at first, there are only point-like manifolds with scalars attached to them.
\\
In these contexts, the $(n, -n)$-PNDP manifolds can play an important role.
\\
In fact \textit{$(n,-n)$-PNDP} appears as a point (point-like), because in general from our interpretation, it is a point (positive and negative dimensions hide each other out and and the total dimension equals zero), but in special it is composed by two manifolds, $B$ and $F$ with nonzero dimensions.
\\
\\
Let's begin with the \textit{MOG  without anomaly case}, and try to consider a tPNDP-manifold to describe the point-like version of the spacetime where to induce a Morris-Thorne wormhole or a Schwarzchild black hole. 
\\
Our aim is to try to demonstrate that with this manifold and the Time, we are able to present a new concept of flat spacetime, that it is point-like version: a PNDP-manifold with $dim=n+m=0$ plus Time.
\\
From 5th equation of (1*), i.e., $f \Delta' f'+(n'-1) |\nabla f|^2 + \lambda f^2 =0$, since $R'=0$ and $\lambda=0$, we obtain $|\nabla f|^2=0$, then $f=$constant.
\\
For this situation, we consider a negative real line $I=-e^{2\Phi(r)}t$ (where $t$ is the Time and $\Phi(r)$ is an adjustable function of $r$) and choose $n=3$, so a (3,-3)-PNDP manifold $T=B \times F$ where $B=\mathbb{R} \times \mathbb{R} \times \mathbb{R}$, but for simplicity we will write as $B=\mathbb{R}^3$ and the same for $F$ that we will write as $F_{(-3)}=\mathbb{R}^3+ E$ (bundle of obstruction), such that $dimF=-3$. Then:
\\
\\
(5) $L=\mathbb{R}^3 \times F_{(-3)} \times I$, with $ds^2=-e^{2\Phi(r)}dt+dx^2+dy^2+dz^2-(d\psi^2 + d\phi^2 + d\sigma^2)_{(-3)}$.
\\
In \cite{Mazharimousavi}, the authors, show that a $(2+1)$-dimensional wormhole can be induced from $(3+1)$-dimensional flat space-time.
\\
They considered the $(3 + 1)$-dimensional Minkowski space-time in the cylindrical coordinates. Then, following the authors' analysis, we can write the metric of $B$ in cylindrical coordinates:
\\
(6)   $ds^2=-e^{2\Phi(r)}dt^2+dr^2+dz^2+r^2d\theta^2-(d\psi^2 + d\phi^2 + d\sigma^2)_{(-3)}$.
\\
If from now we work only temporarily on $B$ plus Time and set $z=\xi(r)$, we find:
\\
\\
(7)   $ds^2=-e^{2\Phi(r)}dt^2+(1+\xi '(r)^2)dr^2+r^2d\theta^2$,
\\
(and if we consider $\Phi(r)=0$ and $\xi '(r)=0$, we have $(2+1)$-Minkowski spacetime), so the form of energy density becomes: $\rho=\frac{\xi ' \xi ''}{r(1+ \xi '^2)^2}$.
\\
The relation (7) is comparable with the Morris-Thorne static wormhole:
\\
(8)   $ds^2=-e^{2\Phi(r)}dt^2+(\frac{1}{1-\frac{b(r)}{r}})dr^2+r^2d\theta^2$, in our specific case we have: $b(r)=r(\frac{\xi '(r)^2}{1+ \xi '(r)^2})$.
\\
\\
It is well-known that the topological structure of a wormhole includes a throat that connects two asymptotically flat spaces and in order to have and maintain this structure, the geometric flare-out condition is necessary, i.e. the minimum dimension to the throat. In the case of Morris-Thorne type wormhole, this condition is given by the huge surface tension compared to the energy density times the square of the speed of light and this implies that, if $r = r_0$ is the location of the throat, then (*) $b(r_0) = r_0$ and (**) for $r> r_0$, $b'(r)< \frac{b(r)}{ r}$. In the new setting, (*) implies that at the throat $\xi ' = \pm \infty$ and (**) states that $\xi ' \xi '' < 0$ for $r > r_0$. Besides these conditions at the throat, we have $z = \xi(r_0) = 0$. 
\\
By the same token, we can construct a Schwarzchild Black Hole, by setting:
\\
$\xi(r)=\frac{4M}{(\frac{2M}{r-2M})^{1/2}}$, and $\Phi(r)=(1/2) ln(\frac{2M}{r}-1)$.
\\
Then we can conclude that from our (5) we can describe a point-like spacetime in which to induce: $((-1)+1)-$spatial dimension Morris-Thorne wormhole or $((-1)+1)-$spatial dimension Schwarzchild black hole.
\\
Our $(n,-n)$-PNDP manifold model consists of two manifolds with nonzero dimensions (one with $n$-dimension and one with $-n$-dimension, where these two manifolds can be thought as a result of intersection between other manifolds). Then we can consider these two manifolds as contained in a "$p$-dimensional BULK", but their product (which generates the $(n,-n)$-PNDP) will create the point-like manifolds.
\\
Returning to the \textit{Graphene wormhole case}, we can carry out the same analysis done above considering $R'=\widetilde R=0$ and $f=const.$, and setting a $(n, -q)$-PNDP-manifolds with $n=4$ and $m=-q=-3$. In this way we obtain the following metric:
\\
(9)   $ds^2=-e^{2\Phi(r)}dt+dx^2+dy^2+dz^2+dw^2-(d\psi^2 + d\phi^2 + d\sigma^2)_{(-3)}$,
\\
that as we saw at the beginning of the section, it is a 1-dimensional manifold plus Time, which in this case we consider to appear as $\mathbb{R}$.
\\
Therefore proceeding with the analysis made for the previous case, but with an extra dimension to the base-manifold $B$, we can obtain:
\\
(10)  $ds^2=-e^{2\Phi(r)}dt^2+(\frac{1}{1-\frac{b(r)}{r}})dr^2+r^2(d\theta^2+sin^2\theta d\zeta^2)-(d\psi^2 + d\phi^2 + d\sigma^2)_{(-3)}$,
\\
which is $[(3+(-3))+1]$-wormhole, i.e., the point-like spatial version plus Time of the graphene wormhole. 
\\
So we can end this "application-session" by saying that the $(n,-n)$-PNDP manifold could be considered as possible mathematical interpretation of  point-like manifolds. 
\\
\\
The PNDP-manifold is a new type of Einstein warped-product manifold that uses the derived-geometry to introduce a concept of virtual dimension, which, from the speculative point of view, is the dimension with which we perceive our PNDP-manifold, the rest is hidden, invisible.
\\
\\
\\
{\bfseries \centerline{5. Conclusions and Remarks}}
\\
\\
 We have introduced this new kind of Einstein warped product manifold:
\\
($\Pi_{i=1}^{q'}B_i \times \Pi_{i=(q'+1)}^{\widetilde q} B_i) \times_f F$, where the base-manifold contains an Einstein-manifold $\Pi_{i=(q'+1)}^{\widetilde q} B_i$ (or $\widetilde B$), and where $\widetilde B$ has the same dimension and the same Einstein-$\lambda$-constant of the PNDP-manifold. 
\\
For this new kind of manifold we have shown a Ricci-flat example with $f$ constant, and we have shown the existence of a non-Ricci-flat case with non-constant warping function $f$. While from the speculative point of view we have relating its dimension to the usual geometric concept of dimension, by means of desuspensions (special projections) $\pi_{(n-d)}: $ PNDP, where in the case $n-d>0$ it is identified with another Einstein-manifold ($\widetilde \Pi_{(i=q'+1)}^{\widetilde q} B_i$) with equal dimension and where the respective scalar curvatures have the same constant $\lambda$, in the case $n-d=0$ it is identified with a point, and in the case $n-d<0$ it is identified with $(n-d)$-th desuspensions of a point.
\\
\\
To conclude, in addition to the possible application seen as point-like versions of the graphene-wormhole and spacetime, the PNDP-manifold could describe a new nature of time. In fact the negative dimension of the fiber-manifold could represent the "past" time, the positive dimension of the base-manifold could represent the time "future" (for example 1-dimensional time "future"), while the "present" time which is instantaneous, could be represented by a point (dimension zero dim(M)=0). It is well-known that Einstein considered the five dimensional Kaluza-Klein theory in his investigation for a unified field theory, i.e., a unification of gravity and electromagnetism, during the years 1938-1943. But by finding a non-singular particle solution, he discarded this theory as an impossible model for a unified field theory. Since then some other theories have been formulated which are based on Kalazu-Klein's theory [\cite{Duf1},\cite{Duf2}, \cite{Duf3}]. In this paper we gave some ideas concerning a supersymmetric Kaluza-Klein theory in which tPNDP-manifold can be considered as hidden extra 1-time-like dimensions in a theory with more than 10 dimensions. To get rid of massless ghosts we have a 4-dimensional spacetime, a spontaneous compactification of ground state is necessary by assuming that the cosmological constant is zero. Indeed, in Einstein theory the vierbein field of the 4-dimensional spacetime has 10 degrees of freedom. But 5-dimensional Kaluza-Klein theory the metric field has 5 more degrees of freedom. 
\\
Suppose that ground state is $\Sigma_5 \times S_1$ with a radius of Planck length. Moreover, let us assume that the metric field is independent of the fifth coordinate. The provided action creates a theory of $\Sigma_4$ which together with the Brans-Dicke scalar field has 15 degrees of freedom. It is well-known that 11-dimensional supergravity is an important theory among other supergravity theories since eleven is the maximal space-time dimension which does not allow particles with helicity greater than two (see: \cite{Nah}). Eleven dimensional supergravity has been discarded as it does not provide a realistic model in 4 dimensions. But it is considered as description of low-energy dynamics of $M$-theory (see: \cite{Duf3}). It is true that the eleven-dimensional super gravity constructed by Cremmer et al. \cite{Crem}, did not consider the chiral phenomenology since it is eleven-dimensional and eleven is odd. But as Wang Mian \cite{Mia} proposed, one can construct a consistent supergravity theory with an extra hidden time-like dimensions. We believe that it is possible to construct a consistent theory of supergravity with hidden extra time-like dimensions, and these hidden extra time-like dimensions are identified with our PNDP-manifold with metric $ds^2=  dt_1^2+dt_2^2-(dt_3^2+dt_4^2)_{-2 }$, which we denote by $M_{PNDP}$.
\\
Our $M_{PNDP}$ is a hidden extra time-like manifold. In fact the "role" of the negative dimension of the fiber-manifold is to hide the dimension of the base-manifold making our $M_{PNDP}$ as a point-manifold (recall, in general, that $dim((n,-n)-PNDP)$ is $dimB+dimF=n - n=0$).
\\
Now it is possible to have the topology of vaccum as $\Sigma_{12} = \Sigma_{4} \times S_{2}\times S_{2}\times S_{1}\times H^{3}\backslash \Gamma$, where $\Sigma_{4}$ is the 4-spacetime with the condition that the cosmological constant is zero. Moreover $H^{3}\backslash \Gamma$ is a hyperbolic manifold and $\Gamma$ is a discrete isometry which is not acting freely. This implies that there is no Killing vectors in this quotient space. 
\\
This suggests that we have no gauge field which is associated with the hidden extra time-like dimensions, and our gravitational model can be written as: $\Sigma_{12}\times M_{PNDP}$. In this way the Time has 5 dimensions, the usual temporal dimension of space-time and 4 extral "hidden" dimensions. We believe that the theory can successfully describe a super-Yang-Mills field theory coupled with gravitation after spontaneous compactification.

\par \bigskip
\
\\Alexander Pigazzini: Mathematical and Physical Science Foundation, Sidevej 5, 4200 Slagelse, Denmark, E-Mail: pigazzini@topositus.com
\\Cenap Ozel: King Abdulaziz University, Department of Mathematics, 21589 Jeddah KSA, E-Mail: cenap.ozel@gmail.com
\\Patrick Linker: Department of Materials Testing, University of Stuttgart, 70569 Stuttgart, Germany, E-Mail: Mrpatricklinker@gmail.com
\\ Saeid Jafari: Mathematical and Physical Science Foundation, Sidevej 5, 4200 Slagelse, Denmark, E-mail: saeidjafari@topositus.com

\end{document}